\documentclass[a4paper,12pt]{article}
\usepackage{epsfig}
\usepackage{amssymb}
\usepackage{latexsym}
\usepackage{url}
\begin{document}
\newenvironment{proof}[1][Proof]
               {\par \normalfont
                \trivlist
                \item[
                \itshape #1{.}]\ignorespaces
               }{\endtrivlist}
\newtheorem{theorem}{Theorem}[section]
\newtheorem{corollary}[theorem]{Corollary}
\newtheorem{lemma}[theorem]{Lemma}
\newtheorem{proposition}[theorem]{Proposition}
\newtheorem{remark}[theorem]{Remark}
\newtheorem{conjecture}[theorem]{Conjecture}
\newcommand{\ind}{1\hspace{-2.5mm}{1}}
\title{Maximizing the size of the giant}
\author{Tom Britton\thanks{Department of
Mathematics, Stockholm University, SE-106 91 Stockholm, Sweden. {\it
E-mail}: \{tom.britton; ptrapman\}@math.su.se}  and Pieter Trapman\footnotemark[1]
}

\maketitle

\begin{abstract}
We consider two classes of random graphs:\\
$(a)$ Poissonian random graphs in which the $n$ vertices in the graph have i.i.d.\ weights distributed as $X$, where $\mathbb{E}(X) = \mu$. Edges are added according to a product measure and the probability that a vertex of weight $x$ shares and edge with a vertex of weight $y$ is given by $1-e^{-xy/(\mu n)}$.\\
$(b)$ A thinned configuration model in which we create a ground-graph in which the $n$ vertices have i.i.d.\ ground-degrees, distributed as $D$, with $\mathbb{E}(D) = \mu$.
The graph of interest is obtained by deleting edges independently with probability $1-p$.

In both models the fraction of vertices in the largest connected component converges in probability to a constant $1-q$, where $q$ depends on $X$ or $D$ and $p$.

We investigate for which distributions $X$ and $D$ with given $\mu$ and $p$, $1-q$ is maximized. We show that in the class of Poissonian random graphs, $X$ should have all its mass at $0$ and one other real, which can be explicitly determined. For the thinned configuration model $D$ should have all its mass at $0$ and two subsequent positive integers.
\end{abstract}

\section{Introduction}

In this paper we study asymptotic properties of some random graphs as $n$, the number of nodes/vertices, tends to infinity. More specifically, we study the size of the largest connected component, within two classes of random graphs. If this largest connected component is of the same order as the number of nodes, then it is called the giant. We find the random graph that has the  largest giant among all random graphs in the class having a pre-defined mean degree.

We consider two types of networks in this paper. Both types of random graphs are frequently used extensions of the classical Erd{\"o}s-R{\'e}nyi graph \cite{Boll01,Durr06}.
\begin{itemize}
\item \textbf{Poissonian random graphs:} Poissonian random graphs were introduced in \cite{Norr06} and are a main example of inhomogeneous random graphs \cite{Boll07}. Our model is slightly different from the model introduced in \cite{Norr06}, but asymptotically (for $n \to \infty$) the fraction of vertices in the largest connected component of the graph will be the same.

We construct the random graph of $n$ vertices as follows. First we assign independent and identically distributed (i.i.d.) weights to the vertices distributed as the non-negative random variable  $X$, with $\mu_X := \mathbb{E}(X) < \infty$. A pair of vertices with weights $x$ and $y$ share an edge with probability $1-\exp[-xy/(n \mu_X)]$,  independent of other edges in the graph. In \cite{Norr06} this probability is $1-\exp[-xy/L_n]$, where $L_n$ is the sum of the $n$ vertex weights.

We note that creating the graph and after that removing edges independently with probability $1-p$, is asymptotically in distribution the same as immediately creating a Poissonian random graph with weight distribution $pX$. Later we will be interested in properties of thinned versions of Poissonian random graphs, but this observation implies that these fall under the same model and hence need no additional analysis.

\item \textbf{Thinned configuration model:} The configuration model
  \cite{Durr06, Moll98} is obtained by assigning i.i.d.\ numbers (distributed
  as the non-negative integer valued random variable $D$) of half
  edges to the $n$ vertices in the graph. We assume $\mu_D = \mathbb{E}(D)<
  \infty$. If the total number of half-edges is odd, then we add one
  half-edge to the final vertex. Then, we pair the vertices uniformly
  at random. The probability that a specific graph is created, is the
  same for all graphs with a given degree sequence. Parallel edges
  (edges with the same end vertices) and self-loops (an edge which
  connect a vertex to itself), might occur, but they will not
  influence the asymptotic fraction of vertices in the largest connected component (when $\mu_D<\infty$ as we
  have assumed), (\cite{Durr06}). We clean the graph, by removing all self loops and merge all parallel edges. It is easy to check that this will not influence the asymptotic degree distribution for $\mu_D<\infty$.
After this, edges are deleted independently with probability $1-p$.
\end{itemize}

Within both of the two random graph models we identify the distribution, $X$ and $D$ respectively, which maximizes the size of the giant among all distribution having some fixed mean $\mu_X$ and $\mu_D$ and $p$ respectively. The problem of identifying which (random) graph has the maximal giant among \emph{all} graphs with a fixed mean degree $\mu_D$ is less interesting in that the giant can make up the whole population whenever $\mu_D\ge 2$, while for $\mu_D<2$, the maximum is obtained in any graph, which is a tree. Similarly the problem of \emph{minimizing} the size of the giant among Poissonian random graphs and thinned configuration models is achieved by chosing $\mathbb{P}(X=0)$ (resp.\ $\mathbb{P}(D=0)$) arbitrary close to 1, which has the effect that the relative size of the giant goes to zero.

Random graphs are interesting in their own right but also have numerous applications. One such application, which has been the inspiration to many scientists, is that of modeling the spread of an infectious disease in a socially organized human community, where the social structure of the community is described by a random graph \cite{Ande99,Ball09,Brit08,Newm02}. The interpretation of thinning the random network when considering the spread of $SIR$ (Susceptible $\to$ Infectious $\to$ Recovered) epidemics with non-random infectious periods on networks is that transmission will take place (if one of the node gets infected) exactly along those edges that are kept, and transmission between one infected node and a susceptible neighbor is modeled to have probability $p$ and being independent of other transmission links. In epidemic terminology our aim is hence to identify the degree distribution with given mean for which the asymptotic fraction of individuals infected, in case of a major epidemic outbreak, is maximal.

\section{Notation and some basic results}

Throughout, we will use $\mathbb{N}$ for the strictly positive integers and $\mathbb{N}_0 = \mathbb{N} \cup 0$ for the non-negative integers. Unless specified otherwise we will use $D$ for a random variable taking values in $\mathbb{N}_0$ and $X$ for a non-negative real valued random variable. The mean of $X$ is denoted by $\mu_X$. A bar above a random variable denotes the size biased variant of the random variable, i.e.
$$\mathbb{P}(\bar{D}=k)=\frac{k \mathbb{P}(D=k)}{\mu_D},$$
or in case of a general positive random variable
$$\mathbb{P}(\bar{X} \leq x)=\frac{\int_0^x y d\mathbb{P}(X \leq y)}{\mu_X} = \frac{\mathbb{E}[X \ind[X \leq x]]}{\mu_X}.$$
The (probability) generating function of an $\mathbb{N}_0$ valued random variable is defined by $f_D(s) := \mathbb{E}[s^D] = \sum_{k=0}^{\infty}s^k \mathbb{P}(D=k)$, for $s \in [0,1]$.
We sometimes use the notation $\bar{f}_D(s)=f_{\bar{D}-1}(s)$. The smallest root of $s = \bar{f}_D(s)$ is denoted by $z_D$ and $q_D$ is defined by $q_D=f_D(z_D)$.
Some well-known facts about generating functions that we will use (and which are easy to check) are:

\begin{enumerate}
\item $f_D(s)$ is analytic on $(0,1)$ and all derivatives of $f_D(s)$ are non-negative.
\item $f_D(1)=1$.
\item $\frac{d}{ds}f_D(s) = \mu_D \bar{f}_{D}(s)$, in particular $\frac{d}{ds}f_D(s)|_{s=1} = \mu_D$. Or equivalently:
\begin{equation}\label{originaltosizebias}
f_D(s) = 1 - \mu_D \int_{s}^1 \bar{f}_{D}(x) dx
\end{equation}
for $s \in [0,1]$.
\item $z_D$ is the extinction probability of a Galton Watson branching process \cite{Jage75} with offspring distribution $\bar{D}-1$ and one ancestor. $q_D$ is the extinction probability of a branching process for which the number of children of the ancestor is distributed as $D$ and all other individuals have offspring distribution $\bar{D}-1$.
\item $q_D$ and $z_D$ are strictly less than 1 if and only if $\mathbb{E}(\bar{D}-1)>1$.
\end{enumerate}

For a non-negative real valued random variable $X$, the distribution of a mixed Poisson($X$) random variable $D$ is given by $\mathbb{P}(D=k) = \mathbb{E}(\frac{X^k}{k!}e^{-X})$. The generating function of this random variable $D$ is given by $f_D(s) = \mathbb{E}(e^{-(1-s)X})$. Furthermore, $\bar{f}_{D}(s) = \mathbb{E}(e^{-(1-s)\bar{X}})$. So the generating function of $\bar{D}-1$ is given by the generating function of a mixed Poisson distribution based on the size biased variant of $X$. We note that $\mu_D = \mu_X$.

In this paper, we consider undirected simple graphs.
A simple graph is a graph with no parallel edges (two or more edges with the same end-vertices) or self-loops (edges connecting a vertex to itself). The degree of a vertex is the number of edges a vertex is adjacent to.

A connected component in a graph is a set of vertices for which there is a path between every pair of vertices in this set. Let $\mathcal{C}^i(n)$ be  the $i$-th largest connected component (in case of a tie, the order of the tied components is uniform at random). 
The number of vertices in a set $\mathcal{S}$ is denoted by $|\mathcal{S}|$.

We consider two types of random graphs (as defined in the introduction):
\begin{itemize}
\item \textbf{Poissonian random graphs:}
For this model there is, for $n \to \infty$, with high probability, at most one giant component \cite{Boll07}, i.e.,
for every $\epsilon>0$ we have,
$$\lim_{n \to \infty} \mathbb{P}(|\mathcal{C}^2(n)| < \epsilon n) = 1.$$
Let $D$ be mixed Poisson ($X$). The fraction of vertices in the giant component is for large $n$, with high probability, close to $1-q_D$. More precise, for every $\epsilon>0$, we have
$$\lim_{n \to \infty} \mathbb{P}(|n^{-1}|\mathcal{C}^1(n)|- (1-q_D)| < \epsilon) =1.$$

\item \textbf{Thinned configuration model:} 
Let $D$ be the degree distribution of the ground-graph and $p$ the thinning parameter.
For technical reasons (see the remark below) we exclude the model in which both $p=1$ and $\mathbb{P}(D=0) + \mathbb{P}(D=2) =1$.
For the thinned configuration model we also have that the probability that there is more than one giant components converges to 0 as $n \to \infty$.
In this class of random graphs, the generating function of the degree distribution of the thinned graph is given by $g_{D,p}(s) =f_D(1-p+ps)$ and $\bar{g}_{D,p}(s) = \bar{f}_D(1-p+ps)$. The mean degree of a randomly chosen vertex is $p \mu_D$. 

Let $z_{D,p}$ be the smallest root of $s=\bar{g}_{D,p}(s)$ and define $q_{D,p}= g_{D,p}(z_{D,p})$.
The fraction of vertices in the giant component is for large $n$, with high probability close to $1-q_{D,p}$. More precise, for every $\epsilon>0$, we have
$$\lim_{n \to \infty} \mathbb{P}(|n^{-1}|\mathcal{C}^1(n)|- (1-q_{D,p})| < \epsilon) =1.$$

\textbf{Remark:} If $p=1$ and $\mathbb{P}(D=0) + \mathbb{P}(D=2) =1$, then there exists $\epsilon>0$ such that  $\lim_{n \to \infty} \mathbb{P}(|\mathcal{C}^2(n)| > \epsilon n) >0.$
In this model the fraction of the vertices in a component of size at least $k$ converges to $\mathbb{P}(D=2)$ for every $k \in \mathbb{N}$. However, the fraction of vertices in the largest component does not converge to $q_{D,p}= \mathbb{P}(D=2)$. 
\end{itemize}

We will show that, for the Poissonian
random graphs with given $\mu_X$, the limiting size of the giant component
is maximized if all vertices have weight $\mu_X$ whenever $\mu_X\ge \mu_c\approx 1.756$. If $\mu_X<\mu_c$, then $X$ should only have mass on $\mu_c$ and $0$. Again we note that thinning with a factor $p$ is equivalent to replacing $X$ by $pX$.

For the thinned configuration model with given $\mu_D$ and $p$, the maximal giant size is obtained if $D$ has all mass on $0$ and two subsequent positive integers $k$ and $k+1$. We were not able to identify a closed formula for $k$, and the exact mass distribution on the three possible atoms.

\section{Poissonian random graph}

Define $\mu_c$ as the largest root of $2 x = e^{x -1/2}$. The numerical value of $\mu_c$ is approximated by $\mu_c \approx 1.756$.

Let $|\mathcal{C}^{1}_X(n)|$ be the size of the giant in the Poissonian random graph with weight distribution $X$ and $n$ vertices. Again $D$ is mixed Poisson($X$). Furthermore, let $q_D$ and $z_D$ be as before.
Let $\mathcal{D}_{\mu}$ be the collection of mixed Poisson random variables with $\mathbb{E}(D)=\mu$.

\begin{theorem}\label{maintheorem}
Let $D^*$ be mixed Poisson($X^*$) and
$$\mathbb{P}(X^*=\max(\mu,\mu_c)) = 1-\mathbb{P}(X^*=0) = \min(1,\mu/\mu_c),$$
then
$$\min_{D \in \mathcal{D}_{\mu}} q_D = q_{D^*}.$$
\end{theorem}

This theorem may be interpreted as follows. If $\mu < \mu_c$, then $D \in \mathcal{D}_{\mu} $ defined via $$\mathbb{P}(X=\mu_c) = 1-\mathbb{P}(X=0) = \mu/\mu_c$$ leads to the Poissonian graph, for which the fraction of vertices in the giant converges in probability (as $n\to\infty$) to the largest limit. If $\mu \geq \mu_c$, then $\mathbb{P}(X=\mu) = 1$ (i.e.\ the Erd{\"o}s-R{\'e}nyi graph \cite{Boll01} with mean degree $\mu$) leads to the largest giant in this class.

We first show by a series of three lemmas that the maximal giant is obtained for a weight distribution with mass only at $0$ and one other real number. After that we show that in this class, $X^*$ leads to the largest giant component, which will complete the proof.

\begin{lemma}\label{holder}
If $A$ and $B$ are positive real valued random variables, then
$$\mathbb{E}(A^kB) \geq \mathbb{E}(B) \left(\frac{\mathbb{E}(AB)}{\mathbb{E}(B)}\right)^k.$$
\end{lemma}

\noindent\textbf{Proof}: Use H{\"o}lders inequality, $\mathbb{E}[XY] \leq \mathbb{E}[(X^a)^{1/a}] + \mathbb{E}[(Y^b)^{1/b}]$, for non-negative random variables $X$ and $Y$ and $a,b>0$ such that $a^{-1} + b^{-1} =1$. Filling in $a=k$, $X = A B^{1/k}$ and $Y=B^{(k-1)/k}$ gives the desired result. $\hfill \Box$

\begin{lemma}\label{belowcrossing}
Let $X$ be a general non negative random variable and $D$ be mixed Poisson($X$). If for some
$s_* \in (0,1)$, $$f_D(s_*) := \mathbb{E}(e^{-X(1-s_*)}) = e^{-\lambda(1-s_*)},$$
then
$\mathbb{E}(e^{-X(1-s)}) \leq e^{-\lambda(1-s)}$ for $s \in [s_*,1]$.
\end{lemma}

\noindent\textbf{Proof}: If $f_D(s)$ crosses $e^{-\lambda(1-s)}$ from below in $s_*$, then we know that
\begin{equation}\label{derivbound}
\frac{d}{ds}f_D(s)|_{s=s^*} \geq \lambda e^{-\lambda(1-s_*)}.
\end{equation}
Furthermore, for $k \in \mathbb{N}$, we have
$$\frac{d^k}{ds^k} f_D(s)|_{s=s_*} = \mathbb{E}(X^ke^{-X(1-s)})  \geq \lambda^k e^{-\lambda(1-s_*)} = \frac{d^k}{ds^k} e^{-\lambda(1-s)}|_{s=s_*}.$$
Here the inequality follows by (\ref{derivbound}) and Lemma \ref{holder} with $A =X$ and $B= e^{-X(1-s_*)}$. Since $f_D(1) = 1$, and
$f_D(s)$ is analytic on $(0,1)$, a Taylor expansion in $s=s_*$ gives that $f_D(s)$ cannot cross $e^{-\lambda(1-s)}$ from below in $s \in (0,1)$. $\hfill \Box$

Let $D:=D(\lambda)$ be mixed Poisson($X$), where $X= X(\lambda)$ is defined by $\mathbb{P}(X=\lambda) = 1-\mathbb{P}(X=0) = \mu/\lambda$.
Let $$f(s;\lambda) := f_D(s)=  1- \mu/\lambda + (\mu/\lambda) e^{-\lambda(1-s)}.$$
Note that $\bar{f}_D(s) = e^{-\lambda(1-s)}$, is the generating function of a Poisson $\lambda$ distribution.
Let $q(\lambda) := q_{D(\lambda)}$ and $z(\lambda) := z_{D(\lambda)}$.

\begin{lemma}\label{mainlemma}
Let $X$ be a general non negative random variable with mean $\mu_X$ and let $D$ be mixed Poisson($X$).\\
(a) If $z_D > z(\mu_X) = q(\mu_X)$, then $q_D > q(\mu_X)$.\\
(b) If $z_D \leq q(\mu_X)$, then $q_D > q(\lambda)$, where $\lambda = - \log[z_D]/(1-z_D)$.
\end{lemma}

\noindent\textbf{Proof:} (a) follows from, $q_D = e^{-X(1-z_D)}$, Jensen's inequality and the fact that $f_D(s)$ is increasing. (b) follows from Lemma \ref{belowcrossing} and the fact that $f_D(s)$ equals $e^{-\lambda (1-s)}=\bar{f}(s,\lambda)$ in $z_D$ and 1. Using (\ref{originaltosizebias}) completes the proof. $\hfill \Box$

\medskip
We now show that among the distributions $(D(\lambda);\lambda>0)$, the fraction of vertices in the giant component will converge in probability to the largest limit for $D^*$.

\begin{lemma}\label{twopoint}
Let $\mathcal{D}'_{\mu}$ be the class of mixed Poisson random variables, where $D$ is mixed Poisson($X$) and $\mathbb{P}(X=\lambda) = 1-\mathbb{P}(X=0) = \mu/\lambda$, where $\lambda \geq \mu$.
For $D \in \mathcal{D}'_{\mu}$. $q_D$ is minimized for $\lambda = \max(\mu,\mu_c)$.
\end{lemma}

\noindent\textbf{Proof:}
We first note that we might assume that $\lambda>1$, otherwise $q_D= 1$ anyway.
In what follows we need that $z(\lambda)$ is differentiable on $(1,\infty)$. We prove this by analyzing the derivative of its inverse $z^{-1}(x)$, and show that it is non-zero on this domain:
From the definition of $z(\lambda)$, we deduce that $z^{-1}(x) = \frac{- \log[x]}{1-x}$. Then
$$\frac{d}{dx}z^{-1}(x) = -\frac{1- x + x \log[x]}{x (1-x)^2}$$
For $x \in (0,1)$, this derivative is strictly negative and finite, by
$$\frac{d}{dx} [-(1- x + x \log[x])] = -\log[x] > 0,$$ and $\left[1- x + x \log[x]\right]_{x=1} = 0$.
This implies that $z(\lambda)$ is differentiable on the domain where it takes values in $(0,1)$, that is on $(1,\infty)$.

Since $\bar{f}(z(\lambda);\lambda) - z(\lambda) = 0$. We obtain by applying the chain-rule
$$0 = \frac{d}{d\lambda} \left[\bar{f}(z(\lambda);\lambda) - z(\lambda)\right] =  \left[\frac{d}{d\lambda} \bar{f}(s;\lambda) + \frac{d}{ds} \bar{f}(s;\lambda) \frac{d}{d\lambda}z(\lambda) - \frac{d}{d\lambda}z(\lambda)\right]_{s=z(\lambda)}.$$
This gives that
$$\frac{d}{d\lambda}z(\lambda) = \left[\frac{\frac{d}{d\lambda} \bar{f}(s;\lambda)}{1-\frac{d}{ds} \bar{f}(s;\lambda)}\right]_{s=z(\lambda)}.$$

Furthermore, $$\frac{d}{d\lambda}q(\lambda) = \frac{d}{d\lambda}f(z(\lambda),\lambda) = \left[\frac{d}{d\lambda} f(s;\lambda) + \frac{d}{ds} f(s;\lambda) \frac{d}{d\lambda}z(\lambda)\right]_{s=z(\lambda)}.$$
Noting that $\frac{d}{ds} f(s;\lambda)|_{s=z(\lambda)} = \mu \bar{f}(z(\lambda);\lambda) = \mu z(\lambda)$, we get
$$\frac{d}{d\lambda}q(\lambda) = \left[\frac{\mu z(\lambda) \frac{d}{d\lambda} \bar{f}(s;\lambda) }{1-\frac{d}{ds} \bar{f}(s;\lambda)} + \frac{d}{d\lambda} f(s;\lambda) \right]_{s=z(\lambda)}.$$
Equating this derivative to 0 and using $f(s;\lambda)= 1- \mu/\lambda + (\mu/\lambda) e^{-\lambda(1-s)}$ and $\bar{f}(s;\lambda)= e^{-\lambda(1-s)}$ gives:
$$ 0 =  \frac{-\mu [z(\lambda)]^2[1-z(\lambda)]}{1-\lambda z(\lambda)} + \frac{\mu}{\lambda^2} [1-z(\lambda)] - \frac{\mu}{\lambda} [1-z(\lambda)]z(\lambda).$$
The solutions of this equation are $z(\lambda) = 1$ and $z(\lambda) = (2\lambda)^{-1}.$
The first solution is of no use, because if $z(\lambda) = 1$, then $q(\lambda)= 1$ as well.
Filling in the second solution in $z(\lambda) = e^{-\lambda(1-z(\lambda))}$, gives $2 \lambda = e^{\lambda -1/2}$.
Because the root of this equality is strictly larger than 1 and because $\lim_{\lambda \to \infty} q(\lambda) = 1$, $q(\lambda)$ takes its minimum on $(1,\infty)$ in this largest root. The lemma follows by observing that $\lambda \geq \mu$. $\hfill \Box$

\medskip
\noindent\textbf{Proof of Theorem \ref{maintheorem}:}
From Lemma \ref{mainlemma} it follows that for any mixed Poisson distribution $D \in \mathcal{D}_{\mu}$, there is a distribution $D' \in \mathcal{D}'_{\mu}$ such that $q_{D'} \leq q_D$. The theorem now follows from Lemma \ref{twopoint}. $\hfill \Box$

\section{Thinned configuration model}
Let $\mathcal{B}_{\mu}$ be the collection of all $\mathbb{N}_0$ valued random variables $D$ with $\mathbb{E}(D) = \mu$. Let $z_{D,p}$ be the smallest root of the equation $s=\bar{g}_{D,p}(s)$ and $q_{D,p} = g_{D,p}(z_{D,p})$.
\begin{theorem}
Let $q^* = \inf_{D \in \mathcal{B}_\mu} q_{D,p}$, then there exist $k \in \mathbb{N}$ and a degree distribution $D^* \in \mathcal{B}(\mu)$, which satisfies $$\mathbb{P}(D^*=0) + \mathbb{P}(D^*=k) + \mathbb{P}(D^*=k+1) =1,$$ such that $q_{D^*,p}=q^*$.
\end{theorem}

\noindent\textbf{Proof:} Let $q:=q_{D,p}$  and $z:=z_{D,p}$. First we show that we always can find a distribution $D'$ with mass only at 0 and 2 subsequent integers, such that the associated $q'=q_{D',p}$ satisfies $q \geq q'$.

Let $q''=q_{D'',p}$ and $z''=z_{D'',p}$, where the degree distribution $D''$ is defined by $\mathbb{P}(D''=\lfloor \mu \rfloor) = 1-\mu + \lfloor \mu \rfloor = 1- \mathbb{P}(D''=\lfloor \mu \rfloor + 1)$. The generating functions associated with the model before thinning are $f''(s)$ and $\bar{f}''(s)$.
We now distinguish two cases.
\begin{itemize}
\item If $z > z''$, let $k = \lfloor \mu \rfloor$ and $\beta = \mu - k$. We rewrite $D$ as a mixture of two distributions $D_1$ and $D_2$, with $\mathbb{E}(D_1) = k$ and $\mathbb{E}(D_2) = k+1$. So $D$ is $D_1$ with probability $1-\beta$ and $D_2$ with probability $\beta$. Then
\begin{eqnarray*}
q= f(1-p+pz) &  = & \mathbb{E}((1-p+pz)^D)\\
\ & = & (1-\beta)\mathbb{E}((1-p+pz)^{D_1}) + \beta \mathbb{E}((1-p+pz)^{D_2})\\
\ & \geq & (1-\beta) (1-p+pz)^{k} + \beta (1-p+pz)^{k+1}\\
\ & = & f_{D''}(1-p+ p z).
\end{eqnarray*}
where the inequality is obtained by Jensen's inequality. By $z > z''$ and the observation that $f_{D''}(s)$ is increasing we note that $$q \geq f_{D''}(1-p+ p z'') = q''$$ and the statement holds for $z > z''$.

\item If $z \leq z''$, then we may chose an integer $k$ and a random variable $D^{\dagger}$ with mass only at $0$, $k$ and $k+1$, such that $\mathbb{E}(D)= \mu$ and $z_{D^{\dagger},p}=:z^{\dagger} = z$. Let $f^{\dagger}(s):= f_{D^{\dagger}}(s) $ be its generating function and $\bar{f}^{\dagger}(s):= \bar{f}_{D^{\dagger}}(s)$. Let $\alpha := \mathbb{P}(\bar{D}^{\dagger}-1=k)$. As is done in the previous section, we use that by (\ref{originaltosizebias})
\begin{equation}\label{niceint}
f(1-p+p z) = 1- \mu \int_{1-p+p z}^1 \bar{f}(s) ds.
\end{equation}
We note that
\begin{eqnarray*}
\frac{d^2}{ds^2} \frac{\bar{f}(s) - \bar{f}^{\dagger}(s)}{s^{k-1}} & = & \frac{d^2}{ds^2} [(\sum_{m=0}^{\infty} \mathbb{P}(\bar{D}-1=m) s^{m-k+1}) -(1-\alpha) + \alpha s] \\
\ & = & \sum_{m=0}^{\infty} \mathbb{P}(\bar{D}-1=m) (m-k+1)(m-k) s^{m-k-1}.
\end{eqnarray*}
Since $m-k+1$ and $m-k$ cannot have opposite signs, all summands are non-negative.
Now if for some $s \in (0,1)$, it holds that $\bar{f}(s) - \bar{f}^{\dagger}(s) = 0$ and $\frac{d}{ds} ((\bar{f}(s) - \bar{f}^{\dagger}(s))s^{-k+1}) \geq 0$, then it is impossible that $$(\bar{f}(1)-\bar{f}^{\dagger}(1))s^{-k+1} =0,$$ which leads to a contradiction. So $\bar{f}^{\dagger}(s) \geq \bar{f}(s)$ on $s \in [z,1]$ and by (\ref{niceint}) the statement of the theorem follows for $z \leq z''$.
\end{itemize}
$\hfill \Box$

An additional question to consider for the thinned configuration model is the following. Let $\bar{\mathcal{B}}_{c}$ be the class of all $p \in (0,1]$ and distributions $D$ with $\mathbb{E}[D]=c/p$. For which $p$ and distribution of $D$ in $\bar{\mathcal{B}}_{c}$ is the (in probability as $n \to \infty$) limit of $n^{-1}|\mathcal{C}^1(n)|$ maximized? Note that in these models the expected number of edges in the thinned graph is kept constant.

If $c>2$, then $p=1$ and $\mathbb{P}(D=\lfloor c \rfloor + 1) = 1 - \mathbb{P}(D=\lfloor c \rfloor) = c - \lfloor c \rfloor$ 
give $q_{D,p}=0$, so  the (in probability as $n \to \infty$) limit of $n^{-1}|\mathcal{C}^1(n)|$  is maximized.

For $c \leq 2$, we first consider the asymptotic branching process and minimize $q_{D,p}$ over $\bar{\mathcal{B}}_{c}$.
The following heuristic argument gives that $q_{D,p}$ is maximized if $p=1$ and $\mathbb{P}(D=0) + \mathbb{P}(D=2) =1$.
Assume that either  $p=1$ or $\mathbb{P}(D=0) + \mathbb{P}(D=2) =1$ does not hold. Then 
$n^{-1}|\mathcal{C}^1(n)|$ converges in probability to $1-q_{D,p}$. The limit of 
$n^{-1}|\mathcal{C}^1(n)|$ is at most $c/2$, because the number of edges in the thinned graph is roughly $cn/2$ and the number of vertices in the giant component of a graph is at most 1 higher than the number of edges. 
Note that for $p=1$ and $\mathbb{P}(D=0) + \mathbb{P}(D=2) =1$, $q_{D,p}=1-c/2$.
However, if $p=1$ and $\mathbb{P}(D=0) + \mathbb{P}(D=2) =1$, $n^{-1}|\mathcal{C}^1(n)|$ does not converge to $c/2$ in probability. Still, the (in probability) limit of $n^{-1}|\mathcal{C}^1(n)|$ can be made arbitrary close to $c/2$ by chosing $\epsilon>0$ arbitrary small and taking $\mathbb{P}(D=2) = c/2- 3 \epsilon$, $\mathbb{P}(D=3) = 2 \epsilon$ and $\mathbb{P}(D=0) = 1-c/2 + \epsilon$.

Formally we can minimze $q_{D,p}$ by applying the change of variables $s=1-pt$ to the right-hand-side of (\ref{niceint}), to obtain 
$$1-q_{D,p} = c \int_0^{1-z_{D,p}} \bar{f}_D(1-pt) dt.$$ Assume that $z_{D,p}<1$, otherwise there will be no giant component anyway. Furthermore, assume that $\mathbb{P}(\bar{D}=1) <1$.
We note that $\bar{f}_D(1-pt) = 1-t$ has roots in $t=0$ and $t=z_{D,p}$. By convexity of $\bar{f}_D(1-pt)$ in $t$, those are the only roots and $\bar{f}_D(1-pt) < 1-t$, for $0 < t < z_{D,p}$. This implies that 
$$1-q_{D,p} \leq c \int_0^{1-z_{D,p}} (1-t) dt \leq c \int_0^{1} (1-t) dt \leq c/2.$$
So, if $c\leq 2$, $1-q_{D,p} < c/2$. On the other hand, the (in probability) limit of $n^{-1}|\mathcal{C}^1(n)|$ can be taken arbitrary close to $c/2$ by taking $p=1$ and  $\mathbb{P}(D=2) = c/2- 3 \epsilon$, $\mathbb{P}(D=3) = 2 \epsilon$ and $\mathbb{P}(D=0) = 1-c/2 + \epsilon$ for arbitrary small $\epsilon>0$. This answers the question posed.

\section*{Acknowledgements}

Both authors are grateful to Riksbankens Jubileumsfond of the Swedish Central Bank.

\end{document}